\documentclass[a4, 12pt]{amsart}
\usepackage{amssymb}
\usepackage{amstext}
\usepackage{amsmath}
\usepackage{amscd}
\usepackage{latexsym}
\usepackage{amsfonts}
\usepackage[all]{xy}

\theoremstyle{plain}
\newtheorem{Theorem}{Theorem}[section]
\newtheorem{Thm}[Theorem]{Theorem}
\newtheorem{thm}[Theorem]{Theorem}

\newtheorem{prop}[Theorem]{Proposition}
\newtheorem{Prop}[Theorem]{Proposition}

\newtheorem{lem}[Theorem]{Lemma}
\newtheorem{Lem}[Theorem]{Lemma}

\newtheorem{cor}[Theorem]{Corollary}
\newtheorem{Cor}[Theorem]{Corollary}

\newtheorem*{thm*}{Theorem}
\newtheorem*{cor*}{Corollary}

\newtheorem*{setting}{Setting}

\newtheorem*{claim*}{Claim}

\theoremstyle{definition}

\newtheorem{Ex}[Theorem]{Example}

\newtheorem{remark}[Theorem]{Remark}
\newtheorem{conj}[Theorem]{Conjecture}

\theoremstyle{remark}

\newtheorem*{cpf}{{\sl Proof of Claim}}

\newtheorem*{tpf2}{{\sl Proof of Theorem 2.1}}

\newcommand{\Proof}{\begin{proof}}
\newcommand{\Qed}{\end{proof}}

\numberwithin{equation}{Theorem}

\newcommand{\rmc}{\mathrm{c}}

\newcommand{\rme}{\mathrm{e}}

\newcommand{\rmr}{\mathrm{r}}

\newcommand{\rmG}{\mathrm{G}}

\newcommand{\calR}{\mathcal{R}}

\newcommand{\fkm}{\mathfrak{m}}

\def\Ht{\operatorname{ht}}

\def\Spec{\operatorname{Spec}}

\def\GCD{\operatorname{GCD}}

\begin{document}

\setlength{\baselineskip}{18pt}

\title{Quasi-socle ideals in Gorenstein numerical semigroup rings}
\author{Shiro Goto}
\address{Department of Mathematics, School of Science and Technology, Meiji University, 1-1-1 Higashimita, Tama-ku, Kawasaki 214-8571, Japan}
\email{goto@math.meiji.ac.jp}
\author{Satoru Kimura}
\address{Department of Mathematics, School of Science and Technology, Meiji University, 1-1-1 Higashimita, Tama-ku, Kawasaki 214-8571, Japan}
\email{skimura@math.meiji.ac.jp}
\author{Naoyuki Matsuoka}
\address{Department of Mathematics, School of Science and Technology, Meiji University, 1-1-1 Higashimita, Tama-ku, Kawasaki 214-8571, Japan}
\email{matsuoka@math.meiji.ac.jp}
\thanks{
\endgraf
}
\maketitle
\begin{abstract}
Quasi-socle ideals, that  is the ideals $I$ of the form $I= Q : \fkm^q$ in Gorenstein numerical semigroup rings over fields are explored,  where $Q$ is a parameter ideal, and $\fkm$ is the maximal ideal in the base local ring, and $q \geq 1$ is an integer. The problems of when $I$ is integral over $Q$ and of when the associated graded ring $\rmG(I) = \bigoplus_{n \geq 0}I^n/I^{n+1}$ of $I$ is Cohen-Macaulay are studied. The problems are rather wild; examples are given. 
\end{abstract}
\section{Introduction}

This paper aims at a study of the Polini-Ulrich Conjecture \ref{conj} (\cite{PU}) of \textit{one-dimensional} case. We shall explore Gorenstein numerical semigroup rings over fields as the test case.
Before stating our own result, let us explain the reason why we are interested in the conjecture of the special case. See Section 2 for the statement of the main Theorem \ref{Main Theorem} of this paper.

Let $A$ be a Cohen-Macaulay local ring with the maximal ideal $\fkm$ and $d=\dim A > 0$.
Let $Q=(a_1, a_2, \cdots ,a_d)$ be a parameter ideal in $A$ and let $q > 0$ be a positive integer. Then we put $I=Q:\fkm^q$ and refer to those  ideals as quasi-socle ideals in $A$.

The study of socle ideals $Q:\fkm$ dates back to the research of L. Burch \cite{B}, where she explored socle ideals of finite projective dimension and gave a very nice characterization of regular local rings (cf. \cite[Theorem 1.1]{GH}). More recently,
A. Corso and C. Polini \cite{CP1, CP2} showed, with the interaction to  linkage theory of ideals, that if $A$ is a Cohen-Macaulay local ring which is \textit{not} regular, one has the equality $I^2=QI$ for every parameter ideal $Q$ in $A$, where $I=Q:\fkm$. Subsequently, the first author and H. Sakurai \cite{GSa1, GSa2, GSa3} showed  the equality $I^2 = QI$ could hold true, where $I = Q:\fkm$,  for numerous parameter ideals $Q$ in $A$, even though the base rings $A$ are not necessarily Cohen-Macaulay. However,
 a more important thing is the following. If $J$ is an equimultiple Cohen-Macaulay ideal of reduction number one, the associated graded ring $\rmG(J) = \bigoplus_{n \ge 0}J^n/J^{n+1}$ is a Cohen-Macaulay ring and, so is the Rees algebra $\calR(J)  = \bigoplus_{n \ge 0} J^n$, provided $\Ht_A 
J \geq 2$. One also knows the number and degrees of the defining equations of $\calR(J)$, so that one can understand fairly explicitly the process of desingularization of $\Spec A$ along the subscheme $\mathrm{V}(J)$.
This observation motivated the ingenious research of C. Polini and B. Ulrich \cite{PU}, where they posed, among many important results, the following conjecture.

\begin{conj}[\cite{PU}]\label{conj}
\textit{Let $(A, \fkm)$ be a Cohen-Macaulay local ring with $\dim A \geq 2$. Assume that $\dim A \geq 3$ when  $A$ is regular. Let $q \geq 2$ be an integer and  let $Q$ be  a parameter ideal in $A$ such that $Q \subseteq \fkm^q$. Then
 $$Q:\fkm^q \subseteq \fkm^q.$$ }
\end{conj}

This conjecture was recently settled by H.-J. Wang \cite{W}, whose beautiful theorem says:

\begin{thm}[\cite{W}]
Let $(A, \fkm)$ be a Cohen-Macaulay local ring with $d = \dim A \ge 2$. Let $q \ge 1$ be an integer and $Q$  a parameter ideal in $A$.
Assume that $Q \subseteq \fkm^q$ and put $I=Q:\fkm^q$.
Then $$I\subseteq \fkm^q, \ \ \fkm^q I = \fkm^q Q,\ \ ~~\operatorname{and}~~\ \ I^2=QI,$$ provided that $A$ is not regular if $d \ge 2$ and that $q \ge 2$ if $d \ge 3$.
\end{thm}

Added to it, the very recent research of S. Goto, N. Matsuoka, and R. Takahashi \cite{GMT}  reported a different approach to the Polini-Ulrich conjecture and proved the following.

\begin{thm}[\cite{GMT}]
Let $(A, \fkm)$ be a Gorenstein local ring with $d=\dim A>0$ and $\rme^0_\fkm(A) \ge 3$, where $\rme^0_\fkm(A)$ denotes the multiplicity of $A$ with respect to $\fkm$.
Let $Q$ be a parameter ideal in $A$ and put $I=Q:\fkm^2$.
Then $\fkm^2 I = \fkm^2 Q$, $I^3=QI^2$, and $\rmG(I)= \bigoplus_{n \geq 0}I^n/I^{n+1}$ is a Cohen-Macaulay ring, so that $\calR(I)=\bigoplus_{n \geq 0}I^n$ is also Cohen-Macaulay, provided $d \ge 3$.
\end{thm}

The researches \cite{W} and \cite{GMT} were independently performed  and their methods of proof are totally different from each other's. Unfortunately, the technique of \cite{GMT} can not go beyond the restrictions that $A$ is a Gorenstein ring, $q=2$, and $\rme_\fkm^0(A) \ge 3$ and however, despite these restrictions, the result \cite[Theorem 1.1]{GMT} holds true    even in the case where $\dim A=1$, while Wang's result says nothing about the case where $\dim A=1$. As is suggested in \cite{GMT}, the one-dimensional case is rather different from higher-dimensional cases and   much more complicated to control.

It seems   natural to ask how one can modify the Polini-Ulrich conjecture, so that it covers also the one-dimensional case.
This question has motivated the present research. We then decided to explore Gorenstein numerical semigroup rings over fields, as the starting point of our investigations, because they are typical one-dimensional Cohen-Macaulay local rings and because higher-dimensional phenomena are often realized, with primitive forms, in those rings of dimension one. We expect, with further investigations,  a generalization of the results  in this paper to higher-dimensional cases and a possible modification of the Polini-Ulrich conjecture, as well.

Let us  explain how this paper is organized.
The statement of the main result Theorem \ref{Main Theorem} and its proof will be found in Section 2.
Theorem 2.1 gives a generalization of \cite[Theorem 1.1]{GMT} in the case where the base rings are numerical semigroup rings. As an application of Theorem 2.1 we will explore in Section 3  numerical semigroup rings $A = k[[t^a, t^{a+1}]]$~$(a > 1)$ over  fields $k$, where $t$ is an indeterminate. We will give a criterion for the ideal $I =(t^s) : \fkm^q$ to be integral over the parameter ideal $(t^s)$ in $A$~(here $q > 0$ is an integer and $0 < s \in H=\left< a, a+1 \right>$, the numerical semigroup generated by $a, a+1$). The problem of when the ring $\rmG (I)$ is Cohen-Macaulay is answered in certain special cases. We agree with the observation  in \cite{GMT} that  the one-dimensional case is wild.
To confirm  this, we will note two examples in Section 4.

\section{The main result and the proof}
Let $0 < a_1 < a_2 < \cdots < a_\ell$ $(\ell \geq 1)$ be integers with $\GCD(a_1, a_2, \cdots ,a_\ell) = 1$.
We put $$H=\left< a_1, a_2, \cdots , a_\ell \right> = \{\sum_{i=1}^{\ell} \alpha_ia_i \mid 0 \le \alpha_i \in \Bbb{Z}\}.$$ Then, because $\GCD(a_1, a_2, \cdots ,a_\ell) = 1$, $H \ni n$ for all $n \in \Bbb Z$ with  $n \gg 0$.
We put $\rmc(H) = \min \{m \in \Bbb{Z} \mid H \ni n \textup{ for all integers } n \ge m\}$, the conductor of $H$.
Let $V=k[[t]]$ be the formal power series ring over a field $k$.
We put $$A=k[[H]] = k[[t^{a_1}, t^{a_2}, \cdots , t^{a_\ell}]] ~~~\subseteq ~~~V.$$ Let $\fkm = (t^{a_1}, t^{a_2}, \cdots , t^{a_\ell})$ be the maximal ideal in $A$. 
Then $A$ is a Cohen-Macaulay local ring with $\operatorname{dim} A = 1$ and $\rme_\fkm^0(A) = a_1$, where $\rme_\fkm^0(A)$ denotes the multiplicity of $A$ with respect to the maximal ideal $\fkm$. The ring $V$ is a module-finite birational extension of $A$. Hence $\overline{A}=V$, where $\overline{A}$ denotes the normalization. We say that the numerical semigroup $H$ is {\it symmetric}, if for every $n \in \Bbb Z$,
$$n \in H ~~\Longleftrightarrow~~\alpha -n \not\in H,$$
where $\alpha  = \rmc (H) - 1$ denotes the Frobenius number of $H$. 
This condition  is equivalent to saying that $A=k[[H]]$ is a Gorenstein ring
 (\cite{HK}).

With this notation we are interested in the problem of when the results of [W] and [GMT] hold true and our result is summarized into the following.

\begin{thm}\label{Main Theorem}
Suppose that $A=k[[H]]$ is a Gorenstein ring. Let $q>0$ be an integer and assume that the following two conditions $(C_1)$ and $(C_2)$ are satisfied for $q$, where $c = \rmc (H)$:
\begin{enumerate}
\item[$(C_1)$] $t^n \in \fkm^q$ for all integers $n \ge c$;
\item[$(C_2)$] Let $n \in H$. Then $n <a_{1}(q-1)$, if $t^n \not\in \fkm ^{q-1}$.
\end{enumerate}
Let $0<s \in H$. Let  $Q=(t^{s})$ and $I=Q:\fkm^q$. Then the following assertions hold true.
\begin{enumerate}
\item[(1)] $\fkm ^q I = \fkm ^{q}Q$ and $Q \cap I^2 = QI$. 
\item[(2)] $I^{2}=QI$, if $s \geq c$.
\item[(3)] $I^{3}=QI^{2}$ and  the associated graded ring $\rmG(I)=\bigoplus_{n \geq 0}I^n/I^{n+1}$ is Cohen-Macaulay, if $s \ge a_{1}(q-1)$.
\end{enumerate}
\end{thm}

Before going ahead, let us note a few remarks on Theorem \ref{Main Theorem}.

\begin{remark} (1) Conditions $(C_1)$ and $(C_2)$ in Theorem \ref{Main Theorem} are naturally satisfied if $a_1 \geq 2$ and $q=1$. We will later show that they are satisfied also in the following two cases.
\begin{enumerate} 
\item[(i)] $A=k[[H]]$ is a Gorenstein ring, $a_{1} \ge 3$, and $q=2$.
\item[(ii)] $\ell=2$, $a_{1}>1$, $a_{2}=a_{1}+1$, and $0<q<a_{1}$. 
\end{enumerate}

(2) In Theorem \ref{Main Theorem} the ring $\rmG (I)$ is not necessarily Cohen-Macaulay and the reduction number $$\rmr_Q(I) = \min \{0 \leq n \in \Bbb Z \mid I^{n+1} = QI^n\}$$ of $I$ with respect to $Q$ can go up, unless $s \geq a_1(q-1)$. See Theorem \ref{no} and Example \ref{ex2}.

(3) Unless condition $(C_2)$ is satisfied, Theorem \ref{Main Theorem} (3) does not hold true in general, although condition $(C_1)$ is satisfied (and hence $I$ is integral over $Q$; cf. Lemma \ref{integral}). See Example \ref{ex1}.
\end{remark}

The rest of this section is devoted to the proof of Theorem \ref{Main Theorem}. Let us restate our setting.
\begin{setting}\label{setting}
Let 
$0 < a_1 < a_2 < \cdots < a_\ell$ $(\ell \geq 1)$ be integers with $\GCD(a_1, a_2, \cdots ,a_\ell) = 1$,

$H = \left<a_i \mid 1 \leq i \leq \ell \right>$, $c = \rmc (H)$, $a = a_1=\min~ [H \setminus \{0\}]$,

$k$ a field, $V=k[[t]]$ the formal power series ring over $k$,

$A=k[[H]] = k[[t^{a_1}, t^{a_2 }, \cdots , t^{a_\ell}]] \subseteq V$, and

$\fkm = (t^{a_1}, t^{a_2}, \cdots , t^{a_\ell})$ the maximal ideal in $A$. 
\end{setting}

We begin with the following. 

\begin{lem}\label{2.3}
Suppose that $a \geq 3$ and let $\alpha \geq a -1$ be an integer. Let $\Lambda=\{n \in \Bbb{Z} \mid 0 \leq n \leq \alpha \}$ and assume that for every $n \in \Lambda$, $n \in H \Leftrightarrow \alpha - n \notin H. $ Then $\alpha = c - 1$, so that $H$ is symmetric.  
\end{lem}
\proof
Let $1 \leq m <a$ be an integer. Then $m \notin H$ and so $\alpha-m \in H$, whence $\alpha+n \in H$ for all $1 \leq n \leq a - 1$. Therefore, since $\alpha \not\in H$, to see that $\alpha = c - 1$, it suffices to show $\alpha +a \in H$. Assume $\alpha + a \not\in H$ and put $\Gamma = \{n \in \Bbb Z \mid 0 \leq n \leq \alpha + a\}$. Let $\Delta = \{ n \in \Bbb Z  \mid \alpha + 1 \leq n \leq \alpha + a -1\}$ and let $$\varphi : \Gamma \cap H \to \Gamma \setminus H, \ \ \ n \mapsto \alpha + a -n.$$ Then, since $\alpha + a \not\in H$ and $\Delta \subseteq H$, we see $$\Gamma \cap H = (\Lambda \cap H) \cup \Delta~~~\operatorname{and} ~~~\Gamma \setminus H = (\Lambda \setminus H) \cup \{\alpha + a\}.$$ Therefore, because  the map $\varphi$ is injective and $\sharp{(\Lambda \cap H)} = \sharp{(\Lambda \setminus H)}$ , we have $$a - 1 = \sharp{\Delta} \leq 1,$$ whence $a \leq 2$, which is impossible. Thus $\alpha + a \in H$ so that $\alpha = c-1$. 
\qed

Let $q>0$ be an integer and let $0<s \in H$. We put $Q=(t^{s})$ and $I=Q:\fkm^q$.
Then $$I= (t^n \mid n \in H, t^n \in I),$$ which is a monomial ideal in $A$. Let $\overline{Q}$ denote the integral closure of $Q$. We then have $\overline{Q} = t^sV \cap A$.

\begin{lem}\label{integral}
Suppose $t^n \in \fkm^q$ for all $n \in \Bbb Z$ such that $n  \ge c$. Then $aq \leq c$ and $I \subseteq \overline{Q}$.
\end{lem}

\proof
We have $aq \leq c$, since $t^{c} \in \fkm^q \subseteq t^{aq}V$ (recall $\fkm \subseteq t^aV$, since $a = \min~ [H \setminus \{0\}$]). Let $n \in H$ and assume $t^n \in I$. We want to show $n \geq s$. Assume the contrary and we see $$(s+c-1)-n= (c-1)+(s-n) \geq c $$ because $s > n$, whence $t^{(s+c-1)-n} \in \fkm^{q}$ by assumption. Therefore, since $t^n \in I = Q: \fkm^q$, we get $$t^{s+c-1}=t^{(s+c-1)-n}t^{n} \in Q=(t^s)$$ whence $t^{c-1}  \in A = k[[H]]$, which is impossible (recall $c = \rmc (H))$. Thus $t^n \in t^sV$, so that $I \subseteq t^sV \cap A = \overline Q$ as is claimed.
\qed

The following result shows that condition $(C_1)$ in Theorem \ref{Main Theorem} is satisfied, if $A=k[[H]]$ is a Gorenstein ring, $a \ge 3$, and $q=2$.

\begin{prop}\label{proposition}
Suppose that $A$ is a Gorenstein ring and let $a \geq 3$. Then $t^n \in \fkm^2$ for all $n \in \Bbb Z$ such that  $n \geq c$. Hence $(t^s) : \fkm^{2} \subseteq \overline{(t^s)}$ for all $0 < s \in H$.
\end{prop}

\proof
We may assume that $H$ is minimally generated by $\{a_i\}_{1 \le i \le \ell}$. Hence $\ell \geq 2$ and $H \ne \left<a_j \mid 1 \leq j \leq \ell, j \ne i \right>$ for all $1 \leq i \leq \ell$. We have $c \geq a \geq 3$, since $0 < c \in H$. Notice that $c > a$. In fact, assume that $c=a$. Then $H \ni n$ for all integers $n \ge a$. Therefore, because $a_i +a_j -a  \ge a$ for all $1 \le i, j \le \ell$, we have $\fkm^2 = t^a\fkm$, so that $\ell_A(A/t^aA) \le 2$, since $A$ is a Gorenstein local ring. This is however impossible, because $\ell_A(A/t^aA) = \rme_{\fkm}^0(A) = a \ge 3$, where $\rme_{\fkm}^0(A)$ denotes the multiplicity of $A$ with respect to $\fkm$.  Hence $c > a$.

Let $n \geq c$ be an integer and assume that $t^{n} \notin \fkm^{2}$. Then $n=a_{i}$ for some $1 \leq i \leq \ell$. We have $i > 1$, since $c > a$. Let 
$$K=\left<a_j \mid 1 \leq j \leq \ell, j \ne i \right>.$$ Then $a_i \notin K$. We have $\GCD(a_j \mid 1 \leq j \leq \ell, j \ne i)=1$. In fact, let $1 \leq m <a$ be an integer. Then $m \not\in H$ but $a_i+m \in H$, since $a_i \geq c$. We write $$a_i+m= \sum_{j=1}^\ell c_j a_j$$ with $0 \leq c_j \in \Bbb Z$. Then $c_i = 0$, because $m \not\in H$. Therefore $a_i +m \in K$ for all $1 \leq m < a$. Hence $a_i+1, a_i+2 \in K$, because $a \geq 3$. Thus $\GCD(a_j \mid 1 \leq j \leq \ell, j \ne i)=1$.

We now apply Lemma \ref{2.3} to the numerical semigroup $K$. 
Let $\alpha = c -1$ and let $0 \leq m \le \alpha$ be an integer. Then, since $0 \leq m < c \leq a_{i}$, we have  $ m \in K=\left<a_j \mid 1 \leq j \leq \ell, j \ne i \right>$, once $m \in H$ (recall that $ a=a_1< a_2< \cdots < a_{\ell}$). Suppose now  that $\alpha -m \not\in K$. Then $\alpha - m \not\in H$ as $0 \le \alpha -m \le \alpha$, whence $m \in H$ because the numerical semigroup $H$ is symmetric, so that we have $m \in K$. Conversely, if  $m \in K$, then $m \in H$, whence $\alpha -m \notin H$ so that  $\alpha - m \not\in K$. Consequently, because $\alpha \ge a$ and $a = \min~[K \setminus \{ 0 \}]$,  by Lemma \ref{2.3} we get $\rmc (K) =\alpha + 1 = c$. Hence $a_i \in K$, because $a_i \geq c$. This is impossible. Thus $t^n \in \fkm^2$ for all integers $n \geq c$. The second assertion follows from Lemma \ref{integral}.
\qed

We are now ready to prove Theorem \ref{Main Theorem}.

\begin{tpf2} 
(1) We will show that $\fkm^{q-1}I \subseteq \fkm^qQ : \fkm$. We put $\Lambda = \{(\alpha_1, \alpha_2, \cdots, \alpha_{\ell}) \in {\Bbb Z}^{\ell} \mid \alpha_i \geq 0~\operatorname{for~all}~1 \leq i \leq \ell ~\operatorname{and}~\sum_{i = 1}^{\ell}\alpha_i = q-1\}$. 
Let $\alpha=(\alpha_1, \alpha_2, \cdots, \alpha_{\ell}) \in \Lambda$ and let $n \in H$ such that $t^n \in I$. Let $\varphi = t^{\sum_{i=1}^{\ell}\alpha_i a_i}{\cdot}t^n$. Then $$\varphi \in \fkm^{q-1}I \subseteq Q : \fkm = (t^s) + (t^{s+c-1}),$$
where the equality $Q : \fkm = (t^s) + (t^{s+c-1})$ follows from the fact that  $A$ is a Gorenstein ring (notice that $t^{s+c-1} \not\in Q = (t^s)$ but $t^m{\cdot}t^{s+c-1} = t^s{\cdot}t^{m + c - 1}\in Q = (t^s)$ for every $0 < m \in H$, because $c = \rmc(H)$ is the conductor of $H$). Consequently $\varphi \in (t^s)$ or $\varphi \in (t^{s+c-1})$, since $\varphi$ is a monomial in $t$. Because $t^m{\cdot}t^{s+c-1} = t^{m + (c-1)}{\cdot}t^s \in \fkm^q Q$ for all $0 < m \in H$ (use condition $(C_1)$; notice that $m + (c-1) \geq c$), we have $\fkm{\cdot}t^{s+c-1} \subseteq \fkm^q Q$. Hence $\fkm \varphi \subseteq \fkm^q Q$ if $\varphi \in (t^{s+c-1})$.

Suppose that $\varphi \in (t^s)=Q$ and write $$\sum_{i=1}^{\ell} \alpha_i a_i + n = h + s$$ with $h \in H$.
Then, since $n \geq s$ by Lemma \ref{integral}, we get  
\begin{eqnarray*}
h &=& \sum_{i=1}^\ell \alpha_i a_i + (n-s)\\
&\ge& \sum_{i=1}^\ell \alpha_i a_i \\
&\ge& a{\cdot}\sum_{i=1}^\ell \alpha_i = a(q-1),
\end{eqnarray*}
so that we have $t^h \in \fkm^{q-1}$ by condition $(C_2)$.
Hence $\varphi = t^{\sum_{i=1}^\ell \alpha_i a_i +n} = t^h{\cdot}t^s\in \fkm^{q-1}Q$ and so $\fkm \varphi \subseteq \fkm^q Q$. Thus $\fkm^{q-1} I\subseteq \fkm^qQ : \fkm$, whence $\fkm^q I = \fkm(\fkm^{q-1}I) \subseteq \fkm^q Q$.

Let us show $Q \cap I^{2}=QI$. Since $\fkm^qI=\fkm^qQ$, we have $\fkm^qI^n=\fkm^qQ^n$ for all $n \in \Bbb Z$. Let $x \in Q \cap I^{2}$ and write $x=t^{s}y$  with $y \in A$. Then for all $\alpha \in \fkm^{q}$, we have $$t^s {\cdot} \alpha y = \alpha x  \in \fkm^{q}I^{2} \subseteq Q^{2} = (t^{2s}).$$ Hence $\alpha y \in Q = (t^s)$ so that we have $y \in Q : \fkm^q = I$. Thus $x \in QI$ whence $Q \cap I^{2}=QI$.

(2) It suffices  to show $I^{2} \subseteq Q$. Let $m,n \in H$ such that $t^m, t^n \in I$. Then $m, n  \ge s \ge c$ by Lemma \ref{integral}. We get $m + n -s \in H$, since $m+ n -s=m +(n-s) \geq c$. Therefore $t^{m}t^{n} = t^{m+n-s}t^s\in Q$, whence $I^{2} \subseteq Q$.

(3) We may assume that $I^{2} \ne QI$. Hence $I^2 \not\subseteq Q$,  because $Q \cap I^{2}=QI$. We have $I \subseteq \fkm^{q-1}$ by condition $(C_2)$, since $s \geq a(q-1)$ and $I \subseteq \overline{Q} \subseteq t^sV$. 
Then, since $I^2 \subseteq \fkm^{q-1} I$, we get 
$$Q \subsetneq Q+I^{2} \subseteq Q:\fkm=Q+(t^{s+c-1}).$$ 
Therefore, since $\ell_A([Q : \fkm]/Q) = 1$ (recall that $A$ is a Gorenstein ring), we have 
$$Q+I^2=Q:\fkm=Q+(t^{s+c-1}),$$ whence $t^{s+c-1} \in I^{2}$ because $t^{s+c-1} \not\in Q$.
Consequently $$I^{2}=(Q \cap I^2) + (t^{s+c-1}) = QI + (t^{s+c-1})$$ because $Q \cap I^2 = QI$, whence  $$I^{3}=QI^{2} + I{\cdot}t^{s+c-1}.$$
Let us check that $I{\cdot}t^{s+c-1} \subseteq QI^2$.  Let $n \in H$ and assume that $t^{n} \in I$. We will show that $t^nt^{s+c-1} \in QI^2$. We may assume that $n > s$. Let $h=(n+s+c-1)-2s=(n-s)+(c-1)$. Then $h \in H$ since $h \geq c$. Therefore $$\alpha t^h{\cdot}t^{2s} =\alpha {\cdot}t^{n}t^{s+c-1} \in \fkm^{q}I^{3} \subseteq Q^3 =(t^{3s})$$ for all $\alpha \in \fkm^{q}$ and so $\alpha t^{h} \in Q$.
Consequently, $t^{h} \in Q : \fkm^q = I$, whence $t^{n}t^{s+c-1}=t^{2s}t^{h} \in Q^{2}I \subseteq QI^{2}$.  Thus $I{\cdot}t^{s+c-1} \subseteq QI^2$ so that $I^3=QI^2$.
Since $I^3=QI^2$ and $Q \cap I^{2}=QI$,  we get $Q \cap I^{i+1}=QI^i$ for all $i \in \Bbb Z$, whence $\rmG(I)$ is a Cohen-Macaulay ring. 
\qed
\end{tpf2}

Combining Proposition \ref{proposition} and Theorem \ref{Main Theorem}, we readily get {\cite[Theorem 1.1]{GMT}} in the case where the base rings are numerical semigroup rings. Notice that condition $(C_2)$ is automatically satisfied for $q = 2$.

\begin{cor}[cf. {\cite[Theorem 1.1]{GMT}}]
Suppose that $A=k[[H]]$ is a Gorenstein ring and that $a \geq 3$. 
Let $0<s \in H$ and put $I=Q:\fkm^2$, where $Q=(t^{s})$. Then the following assertions hold true.
\begin{enumerate}
\item[(1)] $\fkm ^2 I = \fkm ^{2}Q$ and $I^{3}=QI^{2}$. 
\item[(2)] $\rmG(I)=\bigoplus_{n \geq 0}I^n/I^{n+1}$ is a Cohen-Macaulay ring. 
\item[(3)] $I^{2}=QI$, if $s \geq c$.
\end{enumerate}
\end{cor}


\section{The case where $H = \left<a, a+1 \right>$}
In this section let $H=\left<a, a+1 \right>$ with $a \geq 2$. Applying Theorem \ref{Main Theorem}, we shall explore the numerical semigroup $H =\left<a, a+1 \right>$. Let  $c = a(a-1)$, that is the conductor of $H$. Similarly as in Section 2, let $k$ be a field and $A=k[[H]]= k[[t^{a}, t^{a+1}]] \subseteq V$, where $V=k[[t]]$ is the formal power series ring over $k$. We denote by $\fkm=(t^a, t^{a+1})$ the maximal ideal in $A$.

Let $0<s \in H$, $Q=(t^{s})$, and $I=Q:\fkm^q$ with  $q > 0$ an integer. We study the problems of when $I$ is integral over $Q$ and of when the associated graded ring $\rmG(I) =\bigoplus_{n \geq 0}I^n/I^{n+1}$ is a Cohen-Macaulay ring.

Let us begin with the following. 

\begin{Lem}\label{3.1}
The following assertions hold true.
\begin{enumerate}
\item[$(1)$] Let $\ell, i \geq 0$ be integers. Then  $a\ell+i \in H$, if $i \leq \ell$. The converse is also ture, if $i < a$. 
\item[$(2)$] $\fkm^{\ell}=(t^{a \ell +i} \mid 0 \leq i \leq \ell)=(t^n \mid n \in H,~n \geq a \ell)$ for all integers $\ell \geq 0$.
\end{enumerate}
\end{Lem}

\Proof
(1)  If $i \leq \ell$, then certainly $a\ell+i=a(\ell-i)+(a+1)i \in H$. Suppose that $a\ell+i \in H$ and $i < a$. We  write $a\ell+i= \alpha a  +\beta (a+1)$ with $0 \leq \alpha, \beta \in \Bbb Z$. Then $\beta = a [\ell - (\alpha + \beta)] + i$ and so, letting  $m = \ell - (\alpha + \beta)$, we see $m \geq 0$, because $\beta \geq 0$ and $i < a$. Hence 
$$\ell \geq \alpha + \beta \geq \beta = am + i \geq i.$$
Thus $i \leq \ell$.

(2) Let $\ell \geq 0$ be an integer. Then since $$a(\ell -i)+(a+1)i = a\ell + i$$ for all $0 \leq i \leq \ell$, we get $$\fkm^{\ell}={(t^a,t^{a+1})}^{\ell}=(t^{a \ell +i} \mid 0 \leq i \leq \ell).  \leqno{(\sharp{})}$$
To see $\fkm^\ell \supseteq (t^n \mid n \in H$, $n \geq a \ell)$, let $n \in H$ such  that $n \geq a \ell$. We write $n=ap+i$ with $p \geq \ell$ and $0 \leq i<a$. Then $p \geq i$ by assertion (1), so that $t^n = t^{ap + i} \in \fkm^p$ by equality $(\sharp{})$. Hence $t^n \in \fkm^{\ell}$, because $p \geq \ell$. Thus $\fkm^\ell=(t^n \mid n \in H, n \geq a \ell).$
\Qed

\begin{prop}\label{a, a+1}
Conditions $(C_1)$ and $(C_2)$ in $\operatorname{Theorem}$ $\ref{Main Theorem}$ are satisfied for $q$ if and only if $q < a$.
\end{prop}

\Proof
Assume that $q < a$ and let $n \geq c$ be an integer. Then $n \geq aq$, since $q < a$ and $c = a(a-1)$. Hence $t^n \in \fkm^q$ by Lemma \ref{3.1} (2). Let $n \in H$ and assume that $t^n \not\in \fkm^{q-1}$. We then have again by Lemma \ref{3.1} (2) that $n < a(q-1)$. Thus conditions $(C_1)$ and $(C_2)$ in Theorem \ref{Main Theorem} are satisfied. See Lemma \ref{integral} for the {\it only~if} part.
\Qed

The question of when $I$ is integral over $Q$ is now answered in the following way.

\begin{Thm}\label{3.2}
The following three conditions are equivalent to each other.
\begin{enumerate}
\item[$(1)$] $I \subseteq \overline{Q}$.
\item[$(2)$] $\fkm^{q}I=\fkm^{q}Q$.
\item[$(3)$] $q<a$.
\end{enumerate}
\end{Thm}

\Proof
(2) $\Rightarrow$ (1) This is clear and well known (\cite{NR}). 

(3) $\Rightarrow$ (2) This follows from Proposition \ref{a, a+1}. See Theorem \ref{Main Theorem}.

(1) $\Rightarrow$ (3) Assume $q \geq a$. We will check that $s - a \not\in H$. Suppose $s-a \in H$ and let $n \in H$ with $n \ge aq$. Then 
$$n-a \geq aq-a \geq a^2-a=c,$$ whence $(n+s-a)-s=n-a \in H$, so that $t^{n}t^{s-a}=t^{(n+s-a)-s}t^{s} \in Q.$ Because $s-a \in H$ and $\fkm^q = (t^n \mid n \in H, n \geq aq)$ by Lemma \ref{3.1} (2), we get  $t^{s-a} \in Q : \fkm^q = I \subseteq \overline{Q} \subseteq t^sV$ by assumption (1),  which is impossible. Thus $s-a \not\in H$ whence $s > a$. We write $s=a \ell +r$ with $\ell \geq 1$ and $0 \leq r<a$. Then $r > \ell - 1$ by Lemma \ref{3.1} (1) since $s - a = a (\ell - 1) + r \notin H$, while $r \leq \ell$ by Lemma \ref{3.1} (1) since $0 \le r < a$ and $s = a \ell + r \in H$. Thus $r=\ell$ so that  $s=(a+1)\ell$. Hence $\ell < a$ because $s-a < c ~(=a(a-1))$.

Let $n \in H$ with $n \ge aq$. Then
$$a\ell+n-s =n-\ell \geq aq-\ell \geq a^{2}-(a-1)=c+1,$$
whence $a\ell+n-s \in H$, so that $t^{n}t^{a\ell}=t^{a\ell +n-s}t^{s} \in Q$ for all $n \in H$ with $n \geq aq$. Thus $t^{a\ell} \in Q : \fkm^q = I$ since $\fkm^q = (t^n \mid n \in H, n \geq aq)$ by Lemma \ref{3.1} (2). Consequently $t^{a\ell} \in \overline{Q} \subseteq t^sV$ by assumption (1), so that $a\ell \geq s=(a+1)\ell$, which is impossible because $\ell \geq 1$. Thus $q < a$ as is claimed.
\Qed

\begin{Cor}\label{3.3}
Assume that $q < a$. Then the following assertions hold true. 
\begin{enumerate}
\item[(1)] $I^{2}=QI$, if $s \geq aq$.
\item[(2)] $I^3=QI^2$ and $\rmG(I)$ is a Cohen-Macaulay ring, if $s \geq a(q-1)$.
\end{enumerate}
\end{Cor} 

\Proof
Since $q < a$, conditions $(C_1)$ and $(C_2)$ in Theorem \ref{Main Theorem} are satisfied (Proposition \ref{a, a+1}). Hence $Q \cap I^2 =QI$ by Theorem \ref{Main Theorem} (1). Therefore, to see assertion (1), it suffices to show that $I^2 \subseteq Q$. Let $n \in H$ with $t^n \in I$. Then, since $t^n \in \overline{Q} \subseteq t^sV$ by Theorem \ref{3.2}, we have $n \geq s \geq aq$, whence $t^n \in \fkm^q$ by Lemma \ref{3.1} (2). Consequently, $ I \subseteq \fkm^q$, so that we have $I^2 \subseteq \fkm^q I \subseteq Q$ as is required. See Theorem \ref{Main Theorem} (3) for assertion (2). 
\Qed

In order to study $\rmr_Q(I)$ and the question of when $\rmG (I)$ is a Cohen-Macaulay ring in the case where $q<a$, thanks to Corollary \ref{3.3}, we may restrict our attention to the case where $s < aq$.
For the rest of this section we assume that $$q < a ~~\operatorname{and}~~s<aq.$$ We write $s=a\ell +r$ with  $1 \leq \ell < q$ and $0 \leq r<a$. Then  $r \leq \ell$ by Lemma \ref{3.1} (1). We put   $$p=(a-1)+(\ell -q),$$ whence  $\ell \leq p < a-1$.

We shall explore whether $\rmG (I)$ is a Cohen-Macaulay ring or not in  certain special cases. For the purpose we need the following.

\begin{Prop}\label{3.4} $I=Q+\fkm^{p+1}+(t^{ap+i} \mid p-\ell +r<i \leq p).$ In particular, $I=Q+\fkm^{p}$ if $r=0$.
\end{Prop}

\Proof We will show that $I=Q+\fkm^{p+1}+(t^{ap+i} \mid r \leq i \leq p)$. Let $n \in H$. Then by Lemma \ref{3.1} (2) we see
$$t^{n} \in I ~~~\Leftrightarrow ~~~aq+i+(n-s) \in H ~~\operatorname{for~all}  ~~0 \leq i \leq q. \leqno{(\sharp{})}$$ 
Let $n \in H$ such that $n \geq a(p+1)$. Then, since $s=a\ell +r$ and $p=(a-1)+(\ell -q)$, we get  $$aq + (n-s) \geq aq + [a(p+1)-s]=c+(a-r) > c,$$ so that $aq+i + (n-s) \in H$ for all $0 \leq i \leq q$. Hence $t^n \in I$ by $(\sharp{})$ for all $n \in H$ with $n \geq a(p+1)$. Consequently $\fkm^{p+1} \subseteq I$ by Lemma \ref{3.1} (2).

Let $r \leq i \leq p$ and put $n=ap+i$. Then $n \in H$ by Lemma \ref{3.1} (1). We get $aq+(n-s) \geq c$ (use $s=a\ell +r$ and $p=(a-1)+(\ell -q)$), so that $t^{n} \in I$ by $(\sharp{})$. Thus
$$I \supseteq Q+ \fkm^{p+1}+(t^{ap+i} \mid r \leq i \leq p).$$

We put  $K = Q+ \fkm^{p+1}+(t^{ap+i} \mid r \leq i \leq p)$. We will show $I \subseteq K$. Let $n \in H$ with $t^{n} \in I$. We write $n=aq_1+r_1$ with  $q_1 \geq 0$ and $0 \leq r_1<a$. 
Hence $r_1 \leq q_1$ by Lemma \ref{3.1} (1). Then it is clear from the above that $t^n \in K$ if $n \ge ap + r$. Let us  consider the case where $n<ap+r$. We will show that $t^{n} \in Q$. 
We have $n \geq s$, because $t^n \in I \subseteq t^sV$ by Theorem \ref{3.2}. Let $n-s=aq_2+r_2$ with $0 \leq q_2$ and $0 \leq r_2<a$. Then, since $s =a \ell + r$, we have $$aq_2 + r_2 = n-s<ap+r-s=a(p-\ell),$$ whence $q_2 < p-\ell$. Thus  $0< p-\ell-q_2=(a-1)-(q+q_2)$ (recall that $p = (a-1)+(\ell - q)$), so that we have  $$q+q_2<a-1.$$
\begin{claim*} $r_2 \leq q_2$. Hence $n-s =aq_2 + r_2\in H$ by $\operatorname{Lemma}~\ref{3.1}~(1)$.
\end{claim*}
\begin{cpf}
We will show that $r_2+q \leq a-1$. Suppose that $r_2+q > a-1$ and let $i=a-1-r_2$. Then $0 \leq i<q$ and so, since $t^n \in I$, by $(\sharp{})$ we get $$a(q+q_2)+(a-1)=a(q+q_2)+(r_2+i)=aq+i+(n-s) \in H.$$  Therefore  $a-1 \leq q+q_2$ by Lemma \ref{3.1} (1), which is impossible. Hence $r_2+q \leq a-1$. By $(\sharp{})$ we then have $a(q+q_2)+(r_2+q)=aq+q+(n-s) \in H$, because $t^{n} \in I$. Thus  $r_2+q \leq q+q_2$ by Lemma \ref{3.1} (1), so that $r_2 \leq q_2$. Hence $n - s= aq_2+r_2 \in H$ by Lemma \ref{3.1} (1).
\end{cpf}

Thanks to Claim we get $t^{n} \in Q$ if $n < ap + r$. Thus  $I \subseteq K$, whence $I=K$. 

We will show that $I=Q+ \fkm^{p+1}+(t^{ap+i} \mid p-\ell +r<i \leq p)$. Let $r \leq i \leq p-\ell+r$. Then $ap+i \in H$ by Lemma \ref{3.1} (1), since $i \leq p$. Because $s = a \ell + r$, we have $ap+i -s = a(p-\ell)+(i-r) \in H$ (cf. Lemma \ref{3.1} (1); recall that $i-r \leq p-\ell$).  Hence $t^{ap+i} \in Q$. Thus $I \subseteq Q+ \fkm^{p+1}+(t^{ap+i} \mid p-\ell +r<i \leq p)$, whence  $$I=Q+\fkm^{p+1}+(t^{ap+i} \mid p-\ell +r<i \leq p).$$

Because  $I=Q+\fkm^{p+1}+(t^{ap+i} \mid r \leq i \leq p)$, the second assertion  follows from Lemma \ref{3.1} (2).
\Qed

Recall that $\rmr_Q(I) = \min~\{0 \leq n \in \Bbb Z \mid I^{n+1} = QI^n \}$ is the reduction number of $I$ with respect to $Q$. For each $\alpha \in \Bbb R$ let $$\lceil \alpha  \rceil = \min ~\{ n \in \Bbb Z \mid \alpha \leq n \}.$$ With this notation we have the following.

\begin{cor}\label{ok}  Assume that $q=a-1$. Then the following assertions hold true.
\begin{enumerate}
\item[$(1)$] Let $r=0$. Then $I = \fkm^{\ell} $ and $\rmG (I)$ is a Cohen-Macaulay ring.
\item[$(2)$] Let $r = \ell$. Then $\rmG (I)$ is a Cohen-Macaulay ring with $\rmr_Q(I) = \lceil \frac{a-1}{\ell +1} \rceil$. 
\end{enumerate}
\end{cor}

\begin{proof}
(1) We have $p=\ell$ and $s = a \ell$. Hence $I=Q+\fkm^{\ell} = \fkm^{\ell}$ by Proposition \ref{3.4}. Therefore $\rmG (I)$ is a Cohen-Macaulay ring, because so is $\rmG (\fkm) = \bigoplus_{n \geq 0}\fkm^n/\fkm^{n+1}$ (recall that $\rmG (\fkm) \cong k[X,Y]/(Y^a)$, where $k[X,Y]$ denotes the polynomial ring).

(2) We have $p = \ell$ and $I = Q + \fkm^{\ell + 1} $.  Let $x = t^a$ and $y = t^{a+1}$. Hence $t^s = y^{\ell}$. Let $n \geq 1$ be an integer. Then, since $\fkm = (x,y)$, we have $$I^n = QI^{n-1} + \fkm^{n(\ell + 1)} = QI^{n-1} + (x^{n(\ell + 1) - i}y^i \mid 0 \leq i \leq n(\ell + 1)).$$ Let $\ell \leq i \leq n(\ell + 1)$ be an integer. Then, since $$[n(\ell + 1) - i] + (i - \ell)= (n-1)(\ell + 1) + 1, $$ we get $x^{n(\ell + 1) - i}y^i = y^{\ell}{\cdot}x^{n(\ell + 1) - i}y^{i-\ell} \in Q{\cdot}\fkm^{(n-1)(\ell + 1)} \subseteq QI^{n-1}$. Hence $$I^n = QI^{n-1} + (x^{n(\ell + 1) - i}y^i \mid 0 \leq i \leq \ell - 1).$$
 Let $0 \leq i \leq \ell - 1$ be an integer and let 
\begin{eqnarray*}
\varphi &=&a[n(\ell + 1) - i] + (a+1)i - s\\ 
&=& a(n-1)(\ell + 1) + (a + i - \ell).
\end{eqnarray*}
Then, since $0 < a + i - \ell < a$ and $(n-1)(\ell + 1) \geq 0$, we get by Lemma \ref{3.1} (1) that $\varphi \in H$ if and only if $(n-1)(\ell + 1)  \geq a + i - \ell$. When this is the case, we have $x^{n(\ell + 1) - i}y^i = t^{s}{\cdot}t^{\varphi} \in QI^{n-1}$, because $\varphi = a[(n-1)(\ell + 1) - (a + i - \ell)] + (a+1)(a+i-\ell )$ and therefore $t^{\varphi} \in \fkm^{(n-1)(\ell + 1)} \subseteq I^{n-1}$. Let $$\Delta = \{0 \leq i \leq \ell - 1 \mid (n-1)(\ell + 1) < a + i -\ell\}.$$ We then have $$I^n = QI^{n-1} + (x^{n(\ell + 1) - i}y^i \mid i \in \Delta)$$
and summarize this observation into the following.

\begin{claim*} For a given integer $n \geq 1$ the following conditions are equivalent to each other.
\begin{enumerate}
\item[$(1)$] $I^n = QI^{n-1}$.
\item[$(2)$] $Q \supseteq I^n$.
\item[$(3)$] $\Delta = \emptyset$.
\item[$(4)$] $n - 1 \geq  \lceil \frac{a-1}{\ell + 1}  \rceil$. 
\end{enumerate}
Hence $\rmr_Q(I) = \lceil \frac{a-1}{\ell + 1}  \rceil$.
\end{claim*}

\begin{cpf}
The implications $(3) \Rightarrow (1) \Rightarrow (2)$ are clear.

 $(2) \Rightarrow (3)$ See the observation above.

$(3) \Leftrightarrow (4)$ $\Delta = \emptyset$ if and only if $\ell - 1 \not\in \Delta$, and the latter condition is equivalent to saying that $(n-1)(\ell + 1) \geq  a-1$, that is $n - 1 \geq \lceil \frac{a-1}{\ell + 1}  \rceil$.  
\end{cpf}

Now we will show that $Q \cap I^n = QI^{n-1}$. We may assume that $\Delta \ne \emptyset$. Because $$Q \cap I^n = QI^{n-1} + [Q\cap (x^{n(\ell + 1) -i}y^i \mid i \in \Delta)]$$ and the ideals considered are all generated by monomials in $t$, it suffices to show that $$Q \cap (x^{n(\ell + 1) -i}y^i) \subseteq QI^{n-1}$$ for all $i \in \Delta$. Let $R = k[[X, Y]]$ be the formal power series ring over $k$ and let us identify $A = R/(X^{a+1} - Y^a).$
 Let $\overline{*}$ denote the image in $R/(X^{a+1} - Y^a)$. Let  $z \in Q \cap (x^{n(\ell + 1) -i}y^i)$ and write
$$z = y^{\ell}\overline{\eta} = x^{n(\ell + 1) -i}y^i\overline{\rho}$$ with $\eta, \rho \in R$. Then
$$Y^{\ell}\eta = X^{n(\ell + 1) -i}Y^i\rho + (X^{a+1} - Y^a)\delta$$
for some $\delta \in R$. Therefore, since $n(\ell + 1) - i < a + 1$ (recall that $i \in \Delta$), we have 
$$Y^{\ell}(\eta + Y^{a - \ell}\delta )= X^{n(\ell + 1) - i}[Y^i\rho + X^{(a+1) - [n(\ell + 1) - i]} \delta],$$
whence  $$\eta + Y^{a - \ell}\delta = X^{n(\ell + 1) - i}\varepsilon ~~~~ \operatorname{and} ~~~~Y^{i}\rho + X^{(a+1) - [n(\ell + 1) - i ]}\delta = Y^{\ell}\varepsilon$$ for some $\varepsilon \in R$. Here notice that $\delta \in (Y^i)$ and we have 
$$z = y^{\ell}\overline{\eta} \in (y^{a+i}, x^{n(\ell + 1) - i}y^{\ell})$$ in $A= R/(X^{a+1}- Y^a)$. Consequently $z \in QI^{n-1}$, because $a + i - \ell > (n-1)(\ell + 1)$ and $n(\ell + 1) - i \geq (n-1)(\ell + 1)$. Thus $$Q \cap (x^{n(\ell + 1) - i}y^i \mid i \in \Delta) \subseteq Q I^{n-1},$$ whence $Q \cap I^n = QI^{n-1}$ and therefore,  $\rmG (I)$ is a Cohen-Macaulay ring with $\rmr_Q(I) = \lceil \frac{a-1}{\ell + 1}  \rceil$.

\end{proof}

When $q = a-1$ and $r < \ell$, we also have the following estimation of the reduction number $\rmr_Q(I)$ of $I$ with respect $Q$.

\begin{prop}\label{3.7}
Assume that $q = a-1$ and $r < \ell$. Then $\rmr_Q(I) \leq a - \ell$.
\end{prop}

\begin{proof}
Since $q = a-1$ and $s = a \ell + r$, we have $p = \ell$ and $I = (t^{a\ell+ i} \mid r \leq i \leq \ell) + \fkm^{\ell + 1}$ by Proposition \ref{3.4}. 
Therefore, because $$t^{s+i}, ~~t^{a (\ell + 1) + j} \in I ~~\operatorname{for~all}~~0 \leq  i \leq \ell - r ~~\operatorname{and}~~ ~~0 \leq j \leq r - 1,$$ multiplying with  the elements $t^s, t^{s+1} \in I$, we have by induction on $n$ that
$$t^{ns+i}, ~t^{a(\ell + 1) + (n - 1)s + j} \in I^n   ~~\operatorname{for~all}~~0 \leq  i \leq \ell + n  - r - 1 ~~\operatorname{and}~~0 \leq j \leq r - 1$$ for every integer $n \geq 1$.  We now take $n = a - \ell$. Then, since 
$$[a (\ell + 1) + (n-1)s] - [ns + (\ell  + n  - r - 1)] = 1$$
and $(\ell + n - r) + r= \ell + n = a,$ we have $$t^m \in I^{a-\ell}~~\operatorname{ for}~~ s(a - \ell) \leq \forall m \leq s(a - \ell) + (a -1),$$ whence $t^m \in I^{a-\ell}$ for all $m \geq (a- \ell)s$ so that  $t^m \in I^{a-\ell + 1}$ for all $m \geq (a- \ell + 1)s$. Consequently, because $I^m \subseteq t^{ms}V$ for all $m \geq 1$ (recall that  $I \subseteq \overline{Q} \subseteq t^sV$; cf. Theorem \ref{3.2}), we get $I^{a-\ell} = t^{(a-\ell)s}V$ and $I^{a-\ell + 1} = t^{(a-\ell + 1)s}V$, whence  $I^{a-\ell + 1} = QI^{a - \ell}$. Thus $\rmr_Q(I) \leq r - \ell$.
\end{proof}

The following two results show that $\rmG(I)$ is not necessarily a Cohen-Macaulay ring, even though $q < a $.

\begin{thm}\label{no}
 Assume that $q = a-1$ and $r = \ell -1$. Then $Q \cap I^{a-\ell} \ne QI^{(a-\ell) - 1},$ provided  that $\ell \geq 2$ and $a \geq \ell +3$. Hence $\rmG (I)$ is not a Cohen-Macaulay ring and $\rmr_Q(I) = a - \ell$.
\end{thm}

\Proof Since $p = \ell$, by Proposition \ref{3.4} (1) we get $I=Q+\fkm^{\ell+1}+(t^{a\ell+\ell})$. It suffices to show $(t^{a\ell +\ell})^{a- \ell} \in Q$ but $(t^{a\ell +\ell})^{a- \ell} \notin QI^{(a-\ell)-1}$.

Since $s = a \ell + (\ell - 1)$ and $c = a (a-1)$, we have 
$$[(a\ell +\ell)(a-\ell)-s]-c=(\ell-1)a^{2}-(\ell^{2}-1)a-(\ell^{2}+\ell -1)>0$$
(recall that $\ell \geq 2$ and $a \geq \ell +3$), so that $(a\ell +\ell)(a-\ell)-s \in H$ whence $(t^{a\ell +\ell})^{a- \ell} \in Q$. 

To show that $(t^{a\ell +\ell})^{a- \ell} \not\in QI^{(a-\ell)-1}$, we  put $\alpha_1 = a\ell + (\ell -1) ~~(=s)$, $\alpha_2 = a\ell + \ell$, and $\alpha_i = a (\ell + 1) + (i-3)$ for $3 \leq i \leq n$, where $n = \ell + 4$. Then $\alpha_i \in H$ for all $1 \leq i \leq n$ by Lemma \ref{3.1} (1) and $0 < \alpha_i < \alpha_{i+1}$ for all $1 \leq i < n$. 
We furthermore have $I = (t^{\alpha_i} \mid 1 \leq i \leq n)$, because $I=Q+\fkm^{\ell+1}+(t^{a\ell+\ell})$. Hence $$I^{(a-\ell)-1} = (t^{\sum_{i=1}^{n}\alpha_i \beta_i} \mid  0 \leq \beta_i \in \Bbb Z, \ \ \sum_{i=1}^{n}\beta_i = a - \ell -1).$$
We put  $\Lambda = \{(\beta_1, \beta_2, \cdots, \beta_n) \in \Bbb Z^{n} \mid 0 \leq \beta_i \in \Bbb Z, \ \ \sum_{i=1}^{n}\beta_i = a - \ell -1\}$. Assume now that $(t^{a\ell +\ell})^{a- \ell} \in QI^{(a-\ell )-1}$.  Let $$\varphi = (a\ell + \ell)(a - \ell) -s ~(= a^{2}\ell-a\ell^{2}-\ell^{2}-\ell+1).$$  We then have $$\varphi = \sum_{i=1}^{n}\alpha_i\beta_i + h$$
for some $\beta = (\beta_1, \beta_2, \cdots, \beta_n) \in \Lambda$ and $h \in H$.  Then $\beta_1 < a - \ell -1$, since $$\varphi - \alpha_1 (a - \ell - 1) = a - \ell \notin H.$$ Because, for each $1 \leq j \leq n$,  
\begin{eqnarray*}
\varphi \geq \sum_{i=1}^{n}\alpha_i\beta_i&=&\sum_{i \neq j}\alpha_i\beta_i+\alpha_j\beta_j \\ &\geq& \sum_{i \neq j}\alpha_1\beta_i+\alpha_j\beta_j \\
&=&\sum_{i=1}^{n}\alpha_1\beta_i+(\alpha_j-\alpha_1)\beta_j\\
&=&\alpha_1(a-\ell-1)+(\alpha_j-\alpha_1)\beta_j,
\end{eqnarray*}
we notice that if $\beta_j \geq 1$ for some $3 \leq j \leq n$, then 
\begin{eqnarray*}
\varphi &\geq& \alpha_1(a-\ell-1)+(\alpha_j-\alpha_1) \\
&\geq& \alpha_1(a-\ell-1)+(\alpha_3-\alpha_1)\\
&=&a^{2}\ell-a\ell^{2}-\ell^{2}-\ell+2\\
&=&\varphi +1,
\end{eqnarray*}
which is absurd. Thus $\beta_j =0$ for all $3\leq j \leq n$. 
Because $\beta_1 + \beta_2 = a -\ell -1$, we have 
\begin{eqnarray*}
\varphi &=& \alpha_1\beta_1+\alpha_2\beta_2 + h\\ &=&a^{2}\ell-a\ell^{2}-\ell^{2}-\ell-\beta_1+ h\\ & =& (\varphi -1) - \beta_1 + h
\end{eqnarray*} 
whence $h=1+\beta_1 \in H$, which is impossible because  $1 \leq 1+\beta_1<a-\ell<a$. This is a required contradiction and  so $(t^{a\ell +\ell})^{a- \ell} \notin QI^{(a-\ell)-1}$. Hence $\rmG (I)$ is not a Cohen-Macaulay ring, because $Q \cap I^{a-\ell} \ne QI^{(a-\ell) -1}$. We get $\rmr_Q(I) = a - \ell$ for the same reason, because $\rmr_Q(I) \leq
 a - \ell$ by Proposition \ref{3.7}. 
\Qed

\begin{thm}\label{3.9}
Assume that $q = a-1$ and $0 < r < \ell$. Let $k = \ell - r$. Then $Q \cap I^3 \ne QI^2$, if $2\ell + 1 \geq a \geq \ell + k + 2$, whence $\rmG (I)$ is not a Cohen-Macaulay ring.
\end{thm}

\begin{proof}
We have $p = \ell$ and $I = Q + \fkm^{\ell+ 1} + (t^{a \ell +i} \mid r < i \leq \ell) = (t^{a \ell +i} \mid r \leq  i \leq \ell)+(t^{a (\ell + 1) +i} \mid 0 \leq i \leq \ell + 1)$. For each $1 \leq i \leq \ell + k +3$, let $\alpha_i = a \ell + r -1 + i$ if $1 \leq i \leq  k + 1$, and $\alpha_i = a (\ell + 1) + i - (k+2)$ if $k + 2 \leq i \leq \ell + k + 3$. Then $0 < \alpha_i < \alpha_{i+1}$ for all $1 \leq i < \ell + k + 3$ and $I = (t^{\alpha_i} \mid 1 \leq i \leq \ell + k + 3)$.

We put $\varphi= (2a\ell + 2\ell + 1) + s$. Then $\varphi \in H$ by Lemma \ref{3.1} (1), because $\varphi = 3a \ell + [3\ell + (1 - k)]$ and $0 < 3 \ell + (1-k) \leq 3 \ell$. We furthermore have $t^{\varphi} \in I^3$, since $\varphi = (a\ell + r + 1) + 2(a+1)\ell$. We get $t^{\varphi} \in Q$ by Lemma \ref{3.1} (1) as well, because $\varphi - s = 2a \ell + 2\ell + 1 = a (2\ell + 1) + (2\ell + 1 -a) \in H$ (recall that $0 \leq 2 \ell + 1 - a$ by our assumption). We now claim the following, which proves $Q \cap I^3 \ne QI^2$.

\begin{claim*}
$t^{\varphi - s} \not\in I^2$.
\end{claim*}

\noindent
{\it Proof.} Assume that $t^{\varphi - s} \in I^2$ and write $\varphi - s= \alpha_i + \alpha_j + h$ with $1 \leq i \leq j \leq \ell +  k + 3$ and $h \in H$. If $j \geq k + 2$, then 
$$\varphi - s = 2a\ell + 2\ell + 1 \geq \alpha_1 + \alpha_{k+2} = 2a\ell + a + \ell - k.$$ Hence $\ell + 1 \geq a - k \geq \ell + 2$, which is impossible.  Thus $j \leq k +1$, so that 
$$\varphi -s = \alpha_i + \alpha_j + h = 2a\ell + i + j + 2r - 2 + h,$$ whence 
\begin{eqnarray*}
h &=& 2\ell + 1 - [i + j + 2(\ell - k) - 2]\\
&=& 2k - (i + j) + 3 \\
&=& (k + 1 -i ) + (k + 1 - j) + 1 > 0.
\end{eqnarray*}
On the other hand, since $a \geq \ell + k + 2$,  we have 
\begin{eqnarray*}
a-h &=& a - [2k - (i + j) + 3]\\
&\ge& \ell + k + 2 - [2k - (i+j) + 3] \\
&=& (\ell - k) + (i + j -1) \geq 2.
\end{eqnarray*}
Thus $0 < h < a$, which is impossible, because $h \in H$ and $a = \min ~[H \setminus \{ 0 \}]$. Hence $t^{\varphi -s } \not\in I^2$. 
\end{proof}

Thanks to Theorem \ref{3.9}, we have the following, where the {\it if} ~part follows from Corollary \ref{ok}.

\begin{cor}\label{3.10}
Assume that $a \geq 5$ and let $a= 2\ell + 1$ with $\ell \geq 2$ an integer. Let $0 \leq r \leq \ell$ and put $s = a \ell + r$. Then $\rmG (I)$ is a Cohen-Macaulay ring if and only if either $r = 0$ or $r = \ell$, where $I = (t^s) : \fkm^{a-1}$.
\end{cor}

\begin{remark} (1) Corollary \ref {ok} and Theorems \ref{no}, \ref{3.9} give only partial answers, in the case where $q < a$ and $s 
< aq$, to the question of when $\rmG (I)$ is Cohen-Macaulay.

(2) Some of the results of this section are possibly  generalized for  numerical semigroups $H = \left<
a, b \right>$ with $0 < a < b$ and $\operatorname{GCD}(a,b)=1$.

We would like to leave further investigations to interested readers. 
\end{remark}


\section{Examples}
In this section we note two  examples $H=\left<10, 13, 16, 17, 19 \right>, \left< 7, 10, 18, 22 \right>$  of symmetric numerical semigroups, for both of which we consider the ideals $I = (t^s) : \fkm^3$ with $0 < s \in H$. In Example \ref{ex2}, the associated graded rings  $\rmG (I)=\bigoplus_{n \geq 0}I^n/I^{n+1}$ are Cohen-Macaulay  except $ s = 16$, while in Example \ref{ex1}, $\rmG (I)$ are {\it not}  Cohen-Macaulay rings for all but finitely many $0 < s \in H$. Thus, even in the case where $q = 3$, the question of when $\rmG (I)$ is Cohen-Macaulay is rather wild.

In Example \ref{ex1}, condition $(C_1)$ is satisfied but condition $(C_2)$ is not. This shows, to control the Cohen-Macaulay property of the associated graded rings $\rmG (I)$ of $I$, we need both conditions $(C_1)$ and $(C_2)$ in  Theorem \ref{Main Theorem}.

\begin{Ex}\label{ex2}
Let $H = \left<10, 13, 16, 17, 19 \right>$ and $A=k[[t^{10}, t^{13}, t^{16}, t^{17}, t^{19}]]$. Then $\rmc(H)=42$ and $A$ is a Gorenstein local ring, which satisfies,  for $q = 3$, conditions $(C_1)$ and $(C_2)$ in Theorem \ref{Main Theorem}. Let $0<s \in H$. We put $Q=(t^s)$ and $I=Q:\fkm^{3}$. 
Then we have $\fkm^3 I = \fkm^3 Q$, whence $I \subseteq \overline{Q}$, and if $s \ge 10{\cdot}(3-1) = 20$, $\rmG(I)= \bigoplus_{n \geq 0}I^n/I^{n+1}$ is a  Cohen-Macaulay ring, thanks  to Theorem \ref{Main Theorem} (3). When $s < 20$, that is the case  $t^s \notin \fkm^2$, we have the following Table 1, where $\rmr_Q(I) = \min \{0 \leq n \in \Bbb Z \mid I^{n+1} = QI^n\}$. Let us explain how to read the table. The table says, for example, that if $s= 10$, the ideal $I = Q : \fkm^3$ is equal to the maximal ideal $\fkm$ and  generated by five monomials $t^{10}, t^{13}, t^{16}, t^{17}, t^{19}$. We have $\rmr_Q(I) = 3$ and $\rmG (I)$ is  a Cohen-Macaulay ring.

\begin{center}
\begin{table}[h]
\caption{$s < 20$}
\begin{tabular}{|c|c|c|c|}
\hline
$s$ & $I$ & $\rmG(I)$ is CM & $\rmr_Q(I)$\\
\hline
$10$ & $(10,13,16,17,19) = \fkm$ & Yes & $3$\\
\hline
$13$ & $(13,16,19,20,27,34)$ & Yes & $3$\\
\hline
$16$ & $(16,19,23,30,32,34,37)$ & No & $5$\\
\hline
$17$ & $(17,20,23,26,29,35,38)$ & Yes & $2$\\
\hline
$19$ & $(19,25,26,32,33,34,37,40)$ & Yes & $2$\\
\hline
\end{tabular}
\end{table}
\end{center}

Hence $\rmG (I)$ is a Cohen-Macaulay ring if and only if $s \ne 16$. When $s=16$, then  $\rmr_Q(I) = 5$. We have $\ell_A(I^2/QI)=2$, $\ell_A(I^3/QI^2)=\ell_A(I^4/QI^3)=\ell_A(I^5/QI^4)=1$, and $Q \cap I^4  \ne QI^3$.

\end{Ex}

\begin{Ex}\label{ex1}

Let $H=\left< 7, 10, 18, 22 \right>$. Then $\rmc(H)=34$. Let $$A=k[[t^{7}, t^{10}, t^{18}, t^{22}]]~\subseteq~k[[t]],$$ where $k[[t]]$ denotes the formal power series ring over a field $k$. Then $A$ is a Gorenstein ring.
Let $0 < s \in H$. We put $Q=(t^s)$ and $I=Q:\fkm^{3}$. 
Then, since  $t^n \in \fkm^3$ for all $n \in \Bbb Z$ such that $n  \ge 34$, $I \subseteq \overline{Q}$, thanks to Lemma \ref{integral}, and we get the following Tables 2, 3. In Table 2 we assume that $t^s \not\in \fkm^3$. We then have $s = 7, 10, 14, 17, 18, 20, 22, 25,$ or $29$. If $s = 7$, then the ideal $I = Q: \fkm^3$ is generated by four monomials $t^7, t^{18}, t^{20}, t^{22}$, and $\rmG (I) $ is a Cohen-Macaulay ring with $\rmr_Q(I) = 2$ and $\fkm^3 I = \fkm^3 Q$. If $s = 10$, then $\rmG (I)$ is not a Cohen-Macaulay ring, $\rmr_Q(I) = 2$, but $\fkm^3 I \ne \fkm^3 Q$. We have $\rmr_Q(I) = 3$, if $s = 18, 25$.

In Table 3 we assume that $t^s \in \fkm^3$. Then the ideal $I$ is generated by the monomials $t^s, t^{s+4}, t^{s+8}, t^{s+10}, t^{s+13}, t^{s+16}, t^{s+19}$, and $\rmG (I)$ is a Cohen-Macaulay ring  if and only if $s = 21$.   We have $\rmr_Q(I) = 2$ but $\fkm^3 I \ne \fkm^3 Q$ always. 

\begin{center}
\begin{table}[h]
\caption{$Q \not\subseteq \fkm^3$} 
\begin{tabular}{|c|c|c|c|c|}
\hline
$s$ & $I$ & $\rmG(I)$ is CM & $\rmr_Q(I)$ & $\fkm^3I=\fkm^3Q$ \\
\hline
$7$ & $(7,18,20,22)$ & Yes & $2$ & Yes \\
\hline
$10$ & $(10,14,18,29)$ & No  & $2$ & No\\
\hline
$14$ & $(14,18,22,27,30)$ & Yes & $2$ & No \\
\hline
$17$ & $(17,21,25,30,36)$ & No & $2$ & No \\
\hline
$18$ & $(18,22,31,34,37)$ & No & $3$ & No \\
\hline
$20$ & $(20,24,28,36,39)$ & No & $2$ & No \\
\hline
$22$ & $(22,30,35,38,41)$ & Yes & $1$ & No \\
\hline
$25$ & $(25,29,38,40,41,44)$ & No & $3$ & No \\
\hline
$29$ & $(29,37,40,42,45,48)$ & Yes & $1$ & No \\
\hline
\end{tabular}
\end{table}

\begin{table}[h]
\caption{$Q \subseteq \fkm^3$}
\begin{tabular}{|c|c|c|c|c|}
\hline
$s$ & $I$ & $\rmG(I)$ is CM & $\rmr_Q(I)$ & $\fkm^3I=\fkm^3Q$ \\
\hline
$21$ & $(s,s+4,s+8,s+10, \qquad$ & Yes & $2$ & No \\
\cline{1-1}\cline{3-5}
otherwise & $\qquad s+13,s+16,s+19)$ & No & $2$ & No\\
\hline
\end{tabular}
\end{table}
\end{center}

\newpage

Hence $\rmG (I)$ is a Cohen-Macaulay ring if and only if $s =7, 14, 21, 22$, and $29 $. 
\end{Ex}


\end{document}